\documentclass{article}
\usepackage[T2A]{fontenc}
\usepackage[cp1251]{inputenc}
\usepackage[english,russian]{babel}
\usepackage[tbtags]{amsmath}
\usepackage{amsfonts,amssymb,mathrsfs,amscd}
\usepackage[hyper]{msb-a}
\JournalName{}
\numberwithin{equation}{section}

\begin{document}




\begin{fulltext}





	\section{Введение}
	В последние годы были получены принципиальные продвижения в теории интегрируемых биллиардов. Системы из этого важного класса динамических систем с ударами, напомним, имеют ``достаточное количество'' первых интегралов. В недавней программной работе А.Т.Фоменко и В.В.Ведюшкиной \cite{1} сделан подробный обзор результатов и открытых задач об интегрируемых биллиардах и топологии слоений на их фазовом пространстве. Такие системы удается успешно исследовать методами теории топологической классификации слоений Лиувилля и их особенностей, разработанной в работах А.Т.Фоменко, его учеников и соавторов и подробно изложенную в монографии А.В. Болсинова и А.Т.Фоменко \cite{2}.
	
	Класс интегрируемых биллиардов недавно был существенно расширен В.В.Ведюш-\\киной \cite{5} благодаря открытию класса биллиардных книжек (многослойных биллиардов). Эти системы задаются на столах-комплексах, склеенных из плоских столов-областей по общим дугам границы. Имеется отображение (проекция) комплекса на плоскость, а ребра комплекса оснащаются перестановками, задающими переход шара с листа на лист.
	
	Если дуги границы всех плоских столов принадлежат одному и тому же семейству квадрик с общими фокусами (возможно, вырожденному: семейству парабол или семейству окружностей и их радиусов), а перестановки на пересекающихся дугах границы коммутируют, то рассмотренные многослойные биллиарды остаются кусочно-гладко интегрируемыми по Лиувиллю. Фазовая топология таких систем, как оказалось, весьма разнообразна: В.В.Ведюшкиной, А.Т.Фоменко и их соавторами были реализованы многие топологические инварианты особенностей и слоений \cite{5}-\cite{8}, возможных в в гладких и вещественно-анали-\\тических интегрируемых системах.
	
	Как оказалось, эти результаты нетривиально дополняют недавние работы \cite{10}-\cite{12}, посвященные гипотезе Биркгофа об интегрируемых биллиардов и доказательству ее различных формулировок. Опуская некоторые детали, основной эффект можно сформулировать так: интегрируемость нескольких классов биллиардов \textit{требует} принадлежности дуг границы стола семейству софокусных квадрик, а не только \textit{следует} из нее. Как следствие, класс плоских интегрируемых биллиардов без потенциала довольно ``беден'' в смысле послойной гомеоморфности их слоений Лиувилля.

	Найденное В.В.Ведюшкиной обобщение оставляет класс интегрируемых биллиардов достаточно разнообразным, не разрушая при этом доказанную в плоском случае связь интегрируемости с софокусными квадриками. Более того, эта конструкция допускает комбинацию с добавлением потенциала, в поле которого интегрируем биллиард на каждом из плоских листов комплекса.
	
	Отметим, что прежде нами, в основном, изучались свойства слоений на неособых 3-мерных поверхностях постоянно энергии в фазовом пространстве. В настоящей работе сделан следующий шаг: показано, что класс 4-мерных невырожденных фокусных особенностей (встречающихся в полуторических интегрируемых системах, во многих системах из приложений: волчке Лагранжа, сферическом маятнике) целиком реализуется многослойными круговыми биллиардами с потенциалом Гука.
	
	Нашу работу мы посвящаем профессору Владимиру Николаевичу Чубарикову в связи с его семидесятилетием.
	\subsection{Гипотеза А.Т.Фоменко о биллиардах и ее доказательство}
	В программной статье \cite{1} центральное место занимает гипотеза об интегрируемых биллиардах, сформулированная А.Т.~Фоменко. Согласно ней, класс интегрируемых биллиардов ``достаточно широк'' в классе всех интегрируемых систем. Эквивалентностью систем является послойная гомеоморфность их слоений Лиувилля \cite{2} на неособых уровнях $Q^3_h$ энергии $H = h$. Классифицирующим инвариантом Фоменко--Цишанга, см. \cite{3}, является граф с числовыми метками, вершинам которого соответствуют невырожденные особенности ранга~1 в $Q^3_h$.
	
	Ряд положений гипотезы Фоменко уже доказан (см. недавний обзор \cite{4}) в работах В.В.Ведюшкиной, А.Т.Фоменко, И.С.Харчевой и В.А.Кибкало: реализованы все невырожденные особенности ранга 1 (см. \cite{5}), числовые метки инварианта Фоменко-Цишанга (см. \cite{6}, \cite{7}) и базы слоений Лиувилля (см. \cite{8}). По остальным пунктам получены значительные продвижения (реализация многих классов гомеоморфности многообразий $Q^3$ и послойной гомеоморфности слоений на них). Некоторые свойства слоений, как оказалось, не являются сами по себе препятствиями к реализации такого слоения биллиардом (например, среди $Q^3_h$ биллиардов найдены $Q^3$ из класса многообразий Вальдхаузена, но не многообразий Зейферта).
	
	В основе этих результатов лежит открытая В.В.Ведюшкиной \cite{5} конструкция ``биллиардной книжки'', т.е. двумерного стола-комплекса, склеенного из плоских столов-комплексов по общим гладким дугам границы. При этом в качестве каждого листа (2-клетки) книжки берется стол \textit{плоского} интегрируемого софокусного или кругового биллиарда. Тем самым, определена ``проекция'': отображение CW-комплекса на плоскость, являющееся изометрией в ограничении на каждый лист CW-комплекса. Переход биллиардного шара с листа на лист после его удара о 1-клетку комплекса (гладкую дугу склейки нескольких листов-столов по их общей граничной дуге) задается перестановкой на множестве склеенных по ней столов. Пусть для каждой вершины (0-клетки) стола-комплекса перестановки на трансверсально пересекающихся в ней дугах коммутируют. Тогда движение частицы непрерывно. Изоэнергетическая поверхность $Q^3_h$ гомеоморфна гладкому многообразию \cite{9}.
	
	Софокусные или круговые столы плоских биллиардов принадлежат к одному из конечного числа классов (с точки зрения топологии слоения на $Q^3_h$), поскольку углы $3\pi/2$ запрещены для непрерывности движения. Согласно недавно доказанной А.А.Глуцюком \cite{10} версии гипотезы Биркгофа, такие и только такие биллиарды будут полиномиально интегрируемыми (на компактном плоском столе, граница которого является кусочно-гладкой класса $C^2$, и хоть одна ее гладкая дуга не прямолинейна). Близкий результат верен и для плоских биллиардов в постоянном магнитном поле, см. \cite{11}: стол должен быть диском или круговым кольцом с общим центром окружностей. Отметим также доказательство В.Ю.Калошиным и А.Соррентино локальной версии гипотезы Биркгофа \cite{12}.
	
	\subsection{Вопрос реализации 4-мерных особенностей биллиардами}
	 С гипотезой Фоменко тесно связан другой вопрос, поставленный в работе \cite{1}: какие 4-мерные особенности топологии слоений Лиувилля интегрируемых систем встречаются в интегрируемых биллиардах (то есть реализуются последними)? Например, чему послойно гомеоморфно слоение в малой 4-мерной окрестности особой точки? Классы эквивалентности слоений в окрестности точки часто называют \textit{локальными особенностями}. Вопрос о послойной гомеоморфности слоений в окрестности всего слоя, содержащего такие точки, также весьма интересен. Уточним, что теперь от окрестности слоя требуется инвариантность при сдвигах вдоль гамильтоновых полей первых интегралов. Классы полученных слоений называют \textit{полулокальными окрестностями}.
	
	Двумя классами особенностей интегрируемых гамильтоновых систем $M^4, \omega$, $\mathcal{F} = (H, F)$, обнаруженными во многих системах из приложений \cite{2}, являются особенности с вырожденными одномерными орбитами (т.е. содержащие на особом слое окружности критических неботтовских точек ранга 1, см. \cite{2}) или с невырожденными точками ранга 0 (положениями равновесия системы, т.е. точками, где $\mathrm{d} H\,=\,\mathrm{d}\,F = \vec{0}$).
	
	Согласно теореме Элиассона \cite{13}, каждая локальная невырожденная особенность ранга 0 послойно гомеоморфна (и даже симплектоморфна) ровно одной из четырех модельных особенностей (в \cite{2} они были названы особенностями типов центр-центр, центр-седло, седло-седло и фокус-фокус). Классификация \textit{полулокальных} нерасщепляемых особенностей была получена Н.Т.Зунгом \cite{14}, см. также \cite{2}. Все точки ранга 0, попавшие на общий слой, обязаны иметь одинаковый локальный тип. Локальная бифуркационная диаграмма $\Sigma$ этих особенностей состоит из дуг двух трансверсально пересекающихся кривых (в случаях центр-центр, центр-седло, седло-седло) или одной изолированной точки (фокус-фокус). Отметим, что вырожденные особенности ранга~1 обычно находятся в~$\mathcal{F}$-прообразах точек касания и возврата кривых~$\Sigma$.

	Реализация невырожденных особенностей ранга 0 биллиардами была начата В.А.Кибкало в работе \cite{15}. Было показано, что локальные и полулокальные невырожденные особенности типов центр-центр, центр-седло и седло-седло действительно встречаются в слоениях Лиувилля софокусных биллиардов с отталкивающим потенциалом Гука на столах-комплексах. Далее было предложено, как реализовать произвольные особенности типа прямого произведения. Затем в работе \cite{16} были построены биллиардные книжки с потенциалом Гука, особенности которых реализует топологический инвариант каждой полулокальной особенности седло-седло с единственной точкой ранга 0 на слое (включая и случай почти прямого произведения) --- круговую молекулу этой особенности.
	
	В основе данного подхода лежит переход от геодезического потока на эллипсоида к биллиарду при устремлении к нулю малой полуоси эллипсоида. Топологические инварианты слоений Лиувилля потока на эллипсоиде с потенциалом Гука были вычислены И.Ф.Кобцевым в работе \cite{17}. Изучение биллиарда с потенциалом Гука на плоских софокусных столах, отличных от эллипса, было выполнено С.Е.Пустовойтовым, см. \cite{18}, \cite{19}.
	
	В настоящей работе произвольная полулокальная 4-мерная особенность типа фокус-фокус реализована с помощью класса интегрируемых круговых биллиардов с отталкивающим потенциалом Гука на столах-комплексах.
	
	Рассмотрим прямое произведение открытого диска $D^2$ на тор $T^2$ и выберем на одном из торов $n$ гомотопных друг другу нестягиваемых окружностей. Перейдем к фактор-пространству, отождествив точки каждой окружности друг с другом. Полученное слоение с одним особым слоем (тор с $n$ перетяжками, гомотопически эквивалентный букету $n$ сфер $S^2$) послойно гомеоморфно особенности типа фокус-фокус сложности $b$. То же самое число $n$ возникает и как топологический инвариант следующего 3-мерного объекта, не содержащего точек ранга 0. Рассмотрим 3-мерную границу 4-мерной полулокальной особенности фокус-фокус. Ее круговая меченая молекула, т.е. инвариант слоения Лиувилля на этой 3-мерной границе (с базой окружность $S^1$ и слоем тор $T^2$), полностью определяется натуральным числом. Это есть числовая метка монодромии, отвечающая за класс сопряженности матрицы монодромии особенности.
	
	\section{Фокусные особенности и круговые биллиарды с потенциалом}
	В данном разделе для произвольной невырожденной фокусной особенности (возможной в интегрируемых гамильтоновых системах с двумя степенями свободы) построен биллиард на столе-комплексе (склеенном из круговых столов) с отталкивающим потенциалом Гука, имеющий особенность, топологически эквивалентную данной.
	
	\subsection{ Биллиард с потенциалом Гука в круге}	Рассмотрим биллиардный стол $D \subset \mathbb{R}^2(x, y)$, являющийся замкнутым кругом радиуса $R = 1$ с центром в точке $(0, 0)$. Фазовое пространство $M^4$ задается отождествлением пар точка-вектор $(x, v_i)$ из $T^{*} D$ для точек $x\ in \partial D$ граничной окружности и векторов скорости $v_1$ и $v_2$, удовлетворяющих закону ``угол падения равен углу отражения''. Матрица симплектической структуры $\Omega$ имеет стандартный вид в координатах $(x, y, \dot x, \dot y)$ для точки $(x, y)$, не лежащей в малой окрестности границы.
	
	Введем на этом столе центральный отталкивающий потенциал Гука $k (x^2 + y^2) / 2$, где $k <0$. Гамильтониан $H$ такого биллиарда на столе $D$ (который также обозначим $D$) имеет вид
	\[H=\frac{\dot{x}^2+\dot{y}^2}{2}+\frac{k}{2}(x^2+y^2),  \qquad k < 0.\]
	Следующая функция $F$ является первым интегралом биллиарда, т.е. \\ $\omega(\mathrm{sgrad}\,H, \mathrm{sgrad}\,F) = 0$:
	\[F = x \dot{y} - y \dot{x}.\]
	Отражение от границы является упругим и сохраняет функцию $F$. Такая система является кусочно-гладко интегрируемой по Лиувиллю согласно определению А.Т.Фоменко \cite{1}.
	
	\textbf{Утверждение 1.} \textit{Образ фазового пространства $\mathcal{F}(M^4)$ биллиарда $D$ при отображении момента $\mathcal{F} = (H, F)$ и бифуркационная диаграмма $\Sigma$ изображены на рис. 1б. Она состоит из параболы и изолированной точки: начала координат. В прообразе этой точки есть критическая точка ранга 0 типа фокус-фокус. В прообразе точки параболы лежит по одной окружности, и соответствующий 3-атом имеет тип $A$. Прообраз любой другой точки, лежащей между ветвями параболы, гомеоморфен двумерному тору $T^2$. }
	
	\textbf{Доказательство.} 1. Каждая точка из $\Sigma$ принадлежит к одному из двух классов. Первыми являются $\mathcal{F}$-образы критических точек $(x, y, \dot x, \dot y)$ отображения $\mathcal{F}$, лежащих на удалении от границы $D$, т.е. при $x^2 + y^2 < 1$. Вторые являются $\mathcal{F}$-образами точек, где фазовое пространство было склеено по отражению от границы.
	
	Вычислим дифференциалы $d H, d F$ и косые градиеты по формуле $\mathrm{sgrad}\,g= \Omega^{-1}\,\mathrm{d} g$.
	\[\mathrm{d}H\,= \, k x \,\mathrm{d} x \, + \, k y \mathrm{d} y \,+\, \dot{x}\, \mathrm{d} \dot{x} \,+\, \dot{y} \,\mathrm{d} \dot{y}, \qquad \mathrm{sgrad}\,H = -k \dot{x}\,\mathrm{d} x \, - \, k \dot{y}\,\mathrm{d} y \,+\, x \, \mathrm{d} \dot{x} \, +\, y \,\mathrm{d} \, \dot{y},\]
	\[\mathrm{d} F \,= \, \dot{y} \mathrm{d} x \, - \, \dot{x} \, \mathrm{d} y \, - \, y \, \mathrm{d} \dot{x} \,+\, x\, \mathrm{d} \dot{y}, \qquad \mathrm{sgrad}\,F \, = \, y \, \mathrm{d} x \, - \, x \,\mathrm{d} y \,+\, \dot{y} \, \mathrm{d} \dot{x} \, - \, \dot{x} \, \mathrm{d} \dot{y}.\]
	В точке $(x, y, \dot{x}, \dot{y})=(0, 0, 0, 0)$ имеем $\mathrm{d} H\, =\, \mathrm{d} F \, = \, \vec{0}$, т.е. это критическая точка ранга 0. Несложно проверить, что при $k <0$ других критических точек у системы нет.
	
	Определим тип локальной особенности в точке $(0, 0, 0, 0)$, следуя см. \cite{2}; т. 1, гл. 1] . Вычислим линеаризации $A_H$ и $A_F$ векторных полей $\mathrm{sgrad}\,H$ и $\mathrm{sgrad}\,F$:
	\[A_H=\omega^{-1}\mathrm{d}^2H=
	\begin{pmatrix}
	0 && 0 && -k && 0\\
	0 && 0 && 0 && -k\\
	1 && 0 && 0 && 0\\
	0 && 1 && 0 && 0
	\end{pmatrix},\]
\[A_F=\omega^{-1}\mathrm{d}^2F=
	\begin{pmatrix}
	0 && -1 && 0 && 0\\
	1 && 0 && 0 && 0\\
	0 && 0 && 0 && -1\\
	0 && 0 && 1 && 0
	\end{pmatrix}.\]
	Рассмотрим линейную комбинацию $\lambda A_H + \mu A_F$ этих операторов. При $\lambda = \mu = 1$ ее собственные значения имеют вид $\pm i \pm \sqrt{-k}$ для каждой возможной пары знаков. Они попарно различны, отличны от нуля и не являются вещественными или чисто мнимыми при $k < 0$. Следовательно, данная критическая точка является невырожденной и имеет тип фокус-фокус при $k < 0$. Ее образом при отображении $\mathcal{F} = (H, F)$ является начало координат.

	2. Для каждой фиксированной пары значений $(h, f)$ отображения момента $\mathcal{F} = (H, F)$ опишем траектории биллиарда на столе $D$ и область возможного движения, т.е. проекцию слоя Лиувилля на $D$. В полярных координатах $(r, \phi)$ имеем
	$$H=\frac{\dot{r}^2+r^2\dot{\phi}^2}{2}+\frac{kr^2}{2}, \qquad \qquad F=r^2\dot{\phi}.$$
	При $f > 0$ частица движется против часовой стрелки, а при $f < 0$ --- по часовой стрелке. В случае $f =0$ имеем два случая. В первом случае $\dot{\phi}=0$, и частица движения по диаметру круга $D$. Во втором случае  $r=0$: частица находится в центре стола $D$, где координаты $r, \phi$ не регулярны.
	
	Вычислим $\dot{r}$ при $r > 0$. Выразив $\dot{\phi}$ из вида функции $F$ и подставив в уравнение для $H$, имеем
	\[\dot{r}^2 \, = \, 2 h - k r^2  - f^2 / r^2 \, \ge \, 0.\]
	Следовательно, область возможного движения задается неравенствами на $r$
	$$- k r^4 + 2 h r^2 - f^2 \geq 0, \qquad \quad r \leq 1.$$
	Дискриминант левой части $4h^2 - 4 k f^2 \ge 0$ всегда, поскольку $k <0$. При этом меньший корень $\left(- h - \sqrt{h^2 - k f^2}\right) / (-k) \le 0$. Следовательно, область возможного движения является пересечением стола $D$, заданного условием $r \leq 1$, и внешности круга $r \geq r_0 := \sqrt{(- h + \sqrt{h^2 - k f^2} ) / (-k)}$. При следующем условии она непуста и является круговым кольцом (или диском $D$ целиком, в случае $f = 0, h \geq 0$), см. рис 1а:
	$$-h + \sqrt{h^2 - k f^2}\leq -k \qquad \Longleftrightarrow \qquad \quad  h \geq (f^2 + k) / 2.$$
	Правое неравенство задает образ отображения момента $\mathcal{F}$. Его границей является парабола $h=(f^2+k)/2$. Диаграмма $\Sigma$ состоит из нее и точки $(0, 0)$ --- образа критической точки фокус-фокус. Прообраз любой другой точки, лежащей между ветвями этой параболы, гомеоморфен тору. Это несложно проверить, следуя анализу биллиарда без потенциала в круге, см. \cite{21}. 
	
	В случае равенства $h = (f^2 + k) / 2$ имеем, что область возможного движения одномерна и является окружностью -- границей биллиарда. Над каждой ее точкой висит один касательный к ней вектор скорости фиксированной длины, причем его направление зависит от знака $f$. Утверждение доказано. $\square$
	\begin{figure}[h]
		\center{\includegraphics[width=100mm]{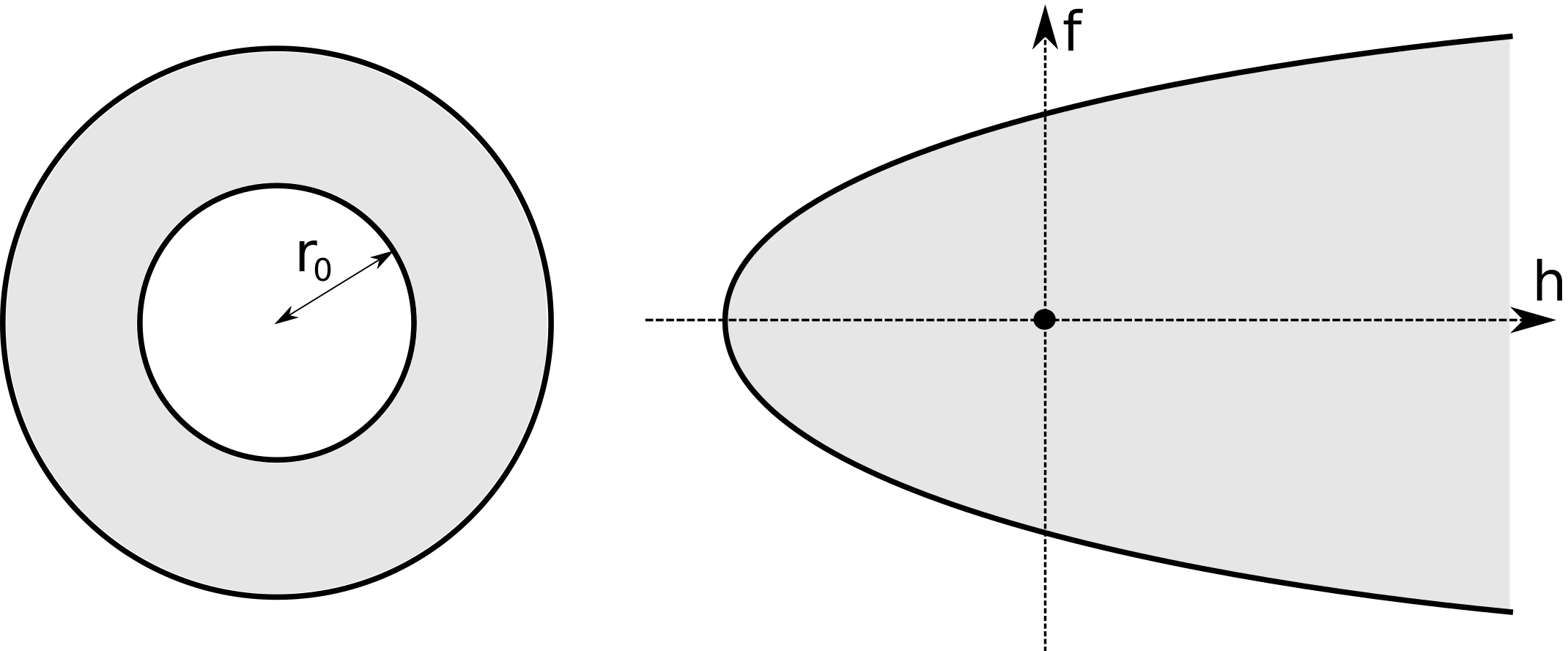}}
		\caption{(а) область возможного движения биллиарда $D$ при фиксированных $h > (f^2 + k) / 2$, (б) образ отображения момента $(H, F)$ для биллиарда $D$ и его бифуркационная диаграмма: парабола $h = (f^2 + k) / 2$ и образ $(0, 0)$ точки фокус-фокус.}
		\label{i1}
	\end{figure}
	
	Заметим, что область возможного движения является всем биллиардным столом $D$ (т. е. $r_0 = 0$), если и только если $h \geq 0$ и $f = 0$. Тогда точка движется вдоль диаметра, пролетая через центр. В случае $h < 0$ и $f = 0$ точка движется вдоль радиуса, не долетая до центра. В случае $f \ne 0$ траектория не лежит на диаметре и является кусочно-гладкой регулярной кривой. При этом направление обхода точки вокруг центра биллиарда (по или против часовой стрелке) определяется знаком $F = f$.
	
	\subsection{ Метка монодромии для биллиардов в круге $\mathbf{D}$ и на комплексах $\Omega_n \mathbf{= n D}$} Для определения полулокального типа особенности в прообразе начала координат у биллиарда $D$ и книжки $\Omega_n$, склеенной из $n$ экземпляров такого стола по циклической перестановке $(1, 2 \dots n)$, мы выполним два независимых вычисления. В данном разделе мы вычислим числовую метку монодромии на круговой молекуле этой особенности, т.е. целочисленный инвариант матрицы перехода между базисами на торе до и после обхода вокруг нее. Для этого мы найдем допустимые базисы и инварианты Фоменко--Цишанга для двух изоэнергетических поверхностей, на которых $H = \pm \varepsilon$ (более точно, вместо левой прямой $H = -\varepsilon$ мы возьмем некоторую параболу на плоскости $Ohf$, на которой $f'_h>0$).
	
	В следующем разделе мы покажем, что инвариантная окрестность особого слоя слоения биллиарда будет послойно гомеоморфна прямому произведению двумерного тора на двумерный открытый диск, в котором один из торов имеет $n$ перетяжек. Иначе говоря, $n$ замкнутых кривых на торе, гомологичных друг другу и не гомологичных нулю в $\pi_1(T^2)$, стянуты в $n$ точек.
	
	\textit{Замечание. } Теоремы 1 и 2 из данного и следующего разделов весьма тесно связаны друг с другом, но используют несколько разные подходы к своим доказательствам. Одной из основных целей наших исследований является расширение набора методов, применяемых для иследования топологии систем биллиардов. Ранее В.В.Ведюшкиной были разработаны эффективные и оригинальные методы изучения топологии, вообще говоря, кусочно-гладких слоений Лиувилля систем биллиардов на 3-уровнях энергии $Q^3_h$. Представляет интерес развитие таких методов в 4-мерном случае. Кроме того, активное расширение класса изучаемых биллиардов (введение систем на столах-комплексах, добавление потенциала или магнитного поля, переход к метрике Минковского, добавление проскальзывания, а также их комбинации, см. \cite{1}) часто дополнительно затрудняет адаптацию и строгое обоснование применимости ряда известных общих методов из теории интегрируемых систем (отметим, например, конструкцию гамильтонова сглаживания В.Ф.Лазуткина, см. \cite{22}).

	\textbf{Теорема 1.} \textit{Числовая метка монодромии особенности в прообразе точки $(0, 0)$ равна $n$ для биллиардной системы с отталкивающим потенциалом Гука на столе $\Omega_n$, склеенном из $n$ экземпляров стола $\Omega_1$ плоского биллиарда внутри круга $D$ по перестановке $(1, 2, \dots, n)$ на их окружности склейки.}
	
	\textbf{Доказательство. } 1. Сначала рассмотрим самый простой случай биллиарда $\Omega_1$ (плоский биллиард). Зафиксируем значение гамильтониана $H = h = - \varepsilon< 0$. На рис. \ref{i2} изображено изменение области возможного движения и поведение траектории с ростом интеграла $F$. Такой перестройке соответствует молекула $A$ --- $A$ с метками $r=\infty$ и $\varepsilon=1$, а соответствующее $Q^3$ гомеоморфно $S^1\times S^2$. Здесь гладкая дуга траектории для $F<0$ стягивается на отрезок радиуса и разворачивается в дугу для $F>0$.
	\begin{figure}[h]
		\center{\includegraphics[width=120mm]{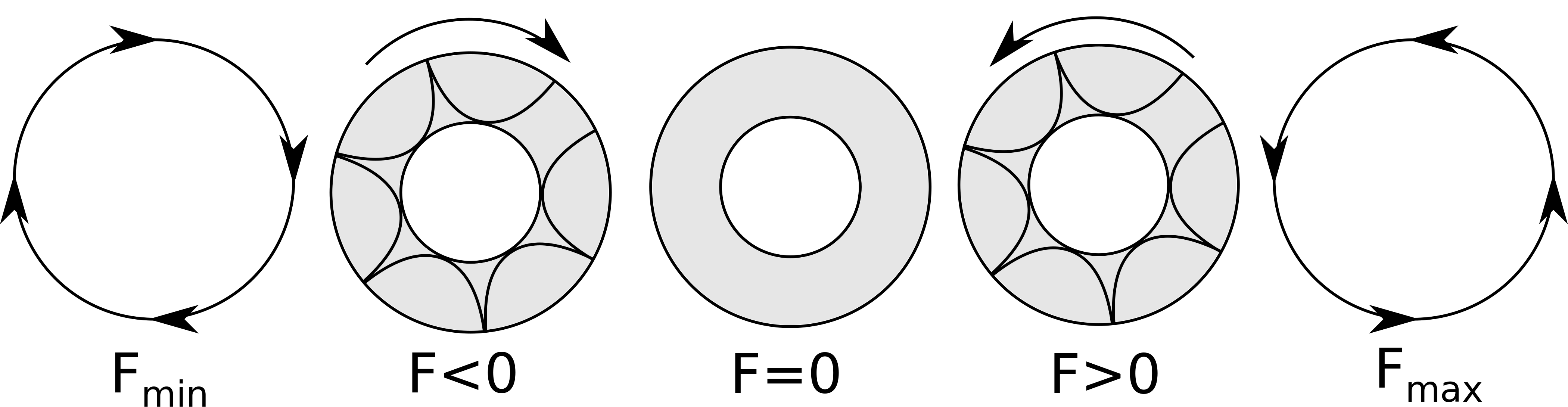}}
		\caption{Область возможного движения и поведение траекторий биллиарда $\Omega_1 = D$ при фиксированном значении $h = -\varepsilon<0$ энергии $H$ в зависимости от значения $f$ интеграла $F$.}
		\label{i2}
	\end{figure}
	
	2. Зафиксируем теперь значение $h = \varepsilon>0$ гамильтониана $H$. В отличие от предыдущего случая, область возможного движения при $f =0$ становится всей биллиардной областью, а гладкие дуги траектории распрямляются в диаметры и затем снова искривляются. Направление обхода точки также определяется знаком $F = f$. Теперь на каждом торе Лиувилля выберем цикл $\gamma$, проекция которого на биллиарде лежит на соседних дугах траектории движения и дуге внутренней границы области возможного движения (голубой цикл на рис. \ref{i3}). Несложно заметить, что при изменении параметра $f$ соответствующий цикл $\gamma$ также меняется непрерывно. Но при $f < 0$ получим $\gamma=\lambda_{-}$ (т.е. цикл, стягивающийся в точку внутри полнотория $A$), а при $f > 0$ получим $\gamma=\lambda_{+}+\mu_{+}$ (сумма исчезающего цикла и дополняющего его до базиса на граничном торе другого полнотория $A$). При этом  для дополняющих циклов $\mu,$ гомологичных внутри полноторий $A$ критическим траекториям (движениям по и против часовой стрелки) выполнено соотношение $\mu_{+}=-\mu_{-}$.  Таким образом, матрица склейки имеет вид $\begin{pmatrix}
	\lambda_{+} \\
	\mu_{+}
	\end{pmatrix}=\begin{pmatrix}
	1 & 1 \\
	0 & -1
	\end{pmatrix}\begin{pmatrix}
	\lambda_{-} \\
	\mu_{-}
	\end{pmatrix}$. Следовательно, такой перестройке соответствует молекула $A$ -- $A$ с метками $r=0$ и $\varepsilon=1$, а соответствующее изоэнергетическое многообразие $Q^3$ гомеоморфно трехмерной сфере $S^3$.
	\begin{figure}[h]
		\center{\includegraphics[width=130mm]{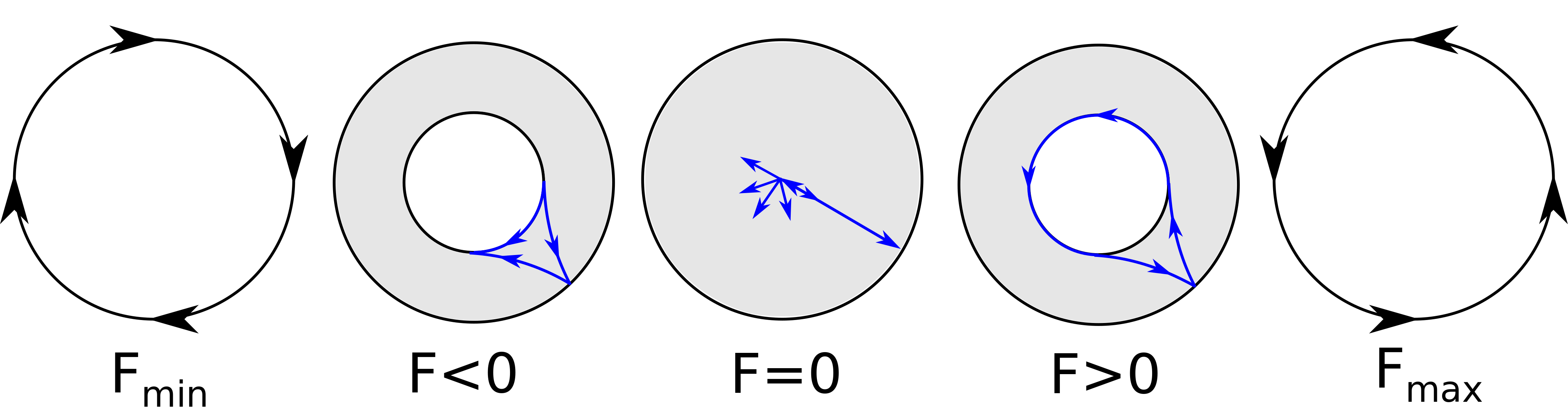}}
		\caption{Область возможного движения и поведение траекторий биллиарда $\Omega_1 = D$ при фиксированном значении $h = \varepsilon> 0$ энергии $H$ в зависимости от значения $f$ интеграла $F$.}
		\label{i3}
	\end{figure}
	
	Покажем, что   полный прообраз точки $(h, f) =(0, 0)$ является ``тором с перетяжкой''. Рассмотрим близкий тор со значениями $h >0, f=0.$ Этот тор отвечает движению вдоль диаметров исходного круга. При стремлении  $h$ к нулю скорость шара вблизи центра круга уменьшается. В пределе $h = 0$ получаем неподвижную точку $(0, 0, 0, 0)$ и некомпактные траектории, направленные к ней или от нее. При $h = 0$ прообразом точки $(0, 0)$ является не окружность (соответсвующая единичным векторам скорости), а одна точка. Стягиваемый в точку при $h \to +0$ цикл был нетривиальным на торе Лиувилля. В результате получен искомый ``тор с перетяжкой''.
	
	3. Теперь сделаем небольшое отступление. Найдем множество значений $h, f$ интегралов $H$ и $F$, для которых радиус внутренней границы равен заданному числу $r_0=c$, где $0 < c< 1$. В этом случае решаем уравнение $\left(-h+\sqrt{h^2-k f^2}\right) / (-k) = c$ и получаем параболу $h = (f^2 + c^2 k) / (2c)$. Теперь в качестве контура обхода вокруг точки фокус-фокус выберем контур, состоящий из дуги этой параболы и отрезка вертикальной прямой (рис. \ref{i4}, а). Изменение траектории движения при обходе этого контура изображено на рис. \ref{i4}, б.
	\begin{figure}[h]
		\center{\includegraphics[width=140mm]{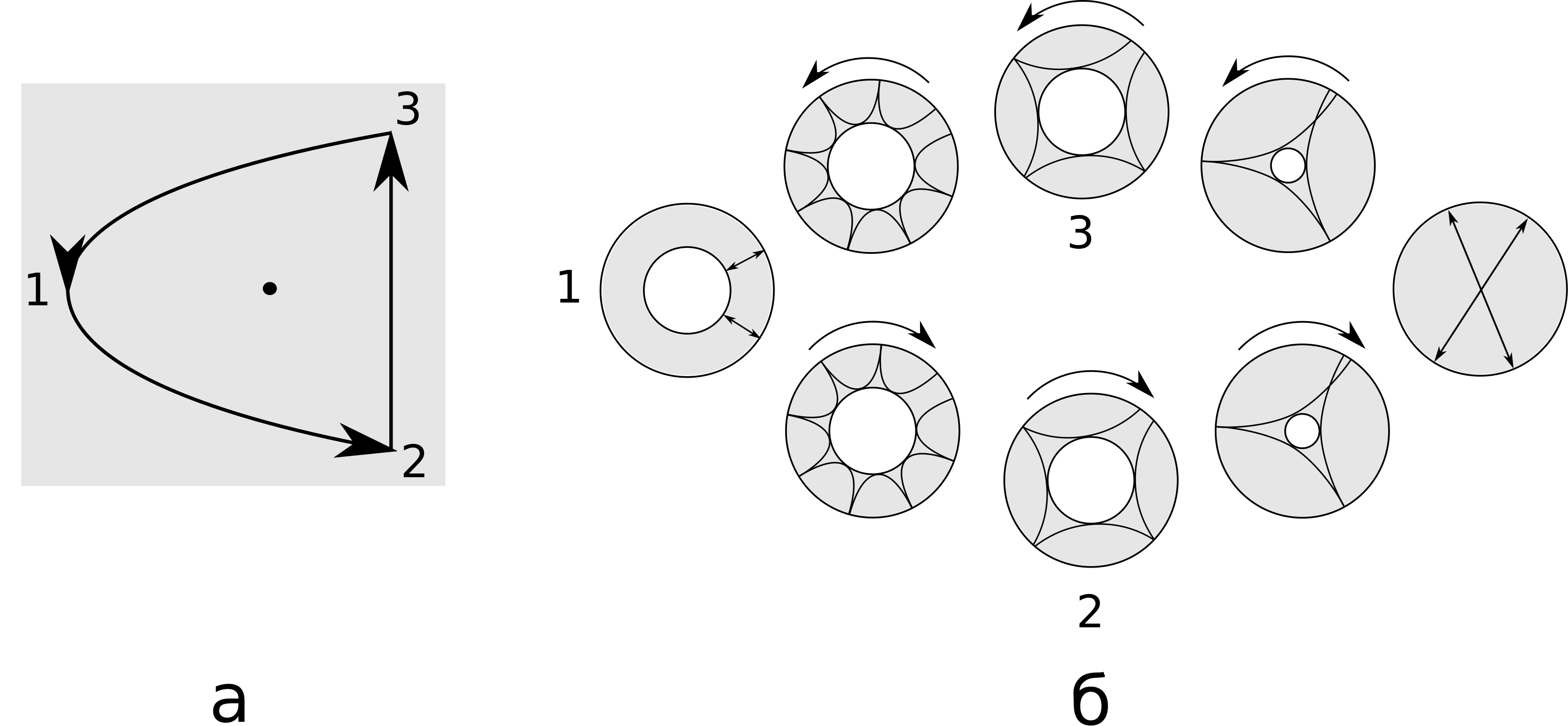}}
		\caption{(а) кривая в $\mathbb{R}^2(h, f)$, обходящая особую точку $\Sigma$, (б) движения частицы в ее прообразе.}
		\label{i4}
	\end{figure}
	На торе Лиувилля, соответствующем точке 1, выберем базис следующим образом: образ цикла $\lambda$ лежит на границе биллиарда, а образ цикла $\mu$ идет по радиусу. Теперь воспользуемся тем же трюком, каким мы пользовались для подсчета матрицы склейки для случая $h > 0$, а именно рассмотрим цикл $\gamma$, выбранный аналогично. В точке 1 $\gamma=\mu$ и остается таким при движении от 1 к 2. При движении от 2 к 3 цикл $\gamma$ изменяется тем же образом, что и раньше, и при переходе от 3 к 1 становится равен $\lambda+\mu$. При этом цикл $\lambda$  при таком же обходе не изменился. Следовательно, матрица монодромии имеет вид $\begin{pmatrix}
	1 & 0 \\
	1 & 1
	\end{pmatrix}$.
	
	4. Мы рассмотрели обычный биллиард в диске с потенциалом. Склеим теперь билллиардную книжку $\Omega_n$ из $n$ копий такого биллиарда по их общей граничной окружности с циклической перестановкой $(1, 2 \dots n)$. Отметим, что в силу того, что область возможного движения всегда ограничена границей биллиарда, такая склейка не добавит системе новых бифуркаций. Следовательно, бифуркационная диаграмма останется прежней.
	
	Тем не менее, прообраз точки $(h, f) = (0, 0)$ теперь устроен сложнее.  Рассмотрим тор Лиувилля, соответствующий значениям $h > 0 ,\ f=0.$ Он склеен из $n$ торов Лиувилля для каждого листа биллиардной книжки, разрезанных и склеенных в единый тор по фиксированным циклам --- прообразам границы кругового биллиарда. Без ограничения общности будем считать эти циклы параллелями торов. Заметим, что эти циклы гомологичны окружностям, проекция которых на биллиард попадает в центр круга. Поэтом при стремлении $h$ к нулю необходимо стянуть каждую из $n$ параллелей в одну точку. В результате получим тор с $n$ "перетяжками".
	
	Как и ранее, значениям $h < 0$ соответствует следующий инвариант слоения Лиувилля на такой изоэнергетической поверхности $Q^3_h$: молекула $A$ -- $A$ с метками $r=\infty$ и $\varepsilon=1$, а $Q^3$ гомеоморфно $S^1\times S^2$. Для $h >0$ зададим цикл $\gamma$ аналогично, но теперь он будет проходить по каждому плоскому листу книжки. Таким образом, при $F>0$ этот цикл равен $\lambda_{+}+n\mu_{+}$ (он обходит вокруг каждого из $n$ центров кругов). Следовательно, этому случаю соответствует молекула $A$ -- $A$ с метками $r=1/n$ и $\varepsilon=1$, а $Q^3$ гомеоморфно $L(1,n)$. Из тех же рассуждений следует, что матрица монодромии при таком же выборе базиса для тора в точке 1 имеет вид $\begin{pmatrix}
	1 & 0 \\
	n & 1
	\end{pmatrix}$. Теорема доказана $\square$
	
	\subsection{ Класс послойной гомеоморфности 4-особенности для биллиарда $\Omega_n$}.
	
	\textbf{Теорема 2.} \textit{Система биллиарда на книжке $\Omega_n$, склеенной из $n$ экземпляров биллиарда в круге $D$ по перестановке $(1, 2, \dots, n)$, имеет в своем слоении Лиувилля полулокальную особенность (содержащую точки $(0, 0, 0, 0)$ из $T^* D$ для каждого из плоских листов), послойно гомеоморфную полулокальной особенности фокус-фокус с $n$ критическими точками на особом слое.}
	
	\textbf{Доказательство.}
	1. Из отсутствия отражений (склеек фазового пространства) в окрестности точки $(0, 0, 0, 0)$ в $T^{*}D$ следует, что система вблизи нее является гладкой. Ранг отображения момента в ней равен нулю, ее невырожденность проверена, и определен тип собственных значений. Иначе говоря, биллиард на столе $D_n$ имеет $n$ локальных особенностей фокус-фокус. Все эти особые точки лежат на одном особом слое, который имеет тип тора с $n$ перетяжками.
	
	2. Остается изучить его 4-мерную окрестность. Она имеет вид $\mathcal{F}^{-1}(V)$ для некоторой малой окрестности $V_\varepsilon \subset \mathbb{R}^2(h, f)$ начала координат $(0, 0)$, гомеоморфной диску. Это верно, поскольку прообраз каждой точки $(h, f)$ связен.
	
	Зафиксируем $y \in [0, 1]$. Для точки $a = (0, y) \in D$ рассмотрим подмножество $A \subset \mathbb{R}^2(\dot x, \dot y)$ точек, образ которых при $\mathcal{F}$ попадает в $V$. Тогда в точке $a$ имеем $2H(a) = k y^2 + \dot{x}^2 + \dot{y}^2$, $F(a) = -y \dot{x}.$
	
	Пусть $y \ne 0$ и $V$ есть прямоугольник $-\varepsilon< h< \varepsilon,  \delta< f < \delta$. Тогда фигура $A$ ограничена, поскольку $\dot{y}^2  = -k y^2 + 2h + f^2/y^2$:
	\[0 \le |\dot{x}| < -\delta/y \qquad  \sqrt{-k y^2 - 2\varepsilon}< |\dot{y}| < \sqrt{-k y^2 + 2\varepsilon + \delta^2/y^2}\]
	Так как $k <0$, то для некоторого $y >0$ фигура $A$ состоит из двух компонент связности, симметричных относительно прямой $\dot{y} = 0$ и имеющих $\dot{y} >0$ или $\dot{y} <0$. При малых $\varepsilon, \delta$ указанный $y$ можно взять меньшим радиуса $1$ круга $D$.
	
	3. Тем самым, пересечение границы $T^*D$ с полулокальной особенностью в полном прообразе $\mathcal{F}^{-1}(V)$ окрестности $V \subset \mathbb{R}^2(h, f)$ имеет две компоненты связности (соответствующие знаку $\dot{y}$, т.е. направлению вектора скорости внутрь круга или наружу). Каждая из них гомеоморфна полноторию $V \times S^1$, где слой $S^1$ над точкой $(h, f) \in V$ есть пересечение слоя Лиувилля $\xi_{h, f}$ с указанной границей.
	
	Склейка двух экземпляров $T^* D$ по паре полноторий (с направлением вектора скорости внутрь круга и наружу соответственно) порождает ориентируемое многообразие: стандартные ориентации этих многообразий согласованы, что легко видеть, например, в склеиваемых точках $x = 0, y = 1, \dot{x}, \pm \dot{y}$.
	
	4. Поскольку $T^* D$ обладает гладкой и симплектической структурой, а также инвариантна относительно вращений, то на $\mathcal{F}^{-1}(V) \subset T^* D$ есть гладкое $S^1$-действие. Оно является гамильтоновым с гамильтонианом $F$. Следовательно, указанное подмножество послойно гомеоморфно элементарному 4-блоку, из $n$ компонент которого склеивается полулокальная особенность типа фокус--фокус, см. \cite{2}, т.1, гл.9. При этом на множество склейки это действие доопределяется по непрерывности из упомянутого выше гамильтонова действия интеграла $F$ на любом из двух склеиваемых по полноторию 4-блоков.
	
	Следовательно, полулокальная особенность биллиарда с отталкивающим потенциалом Гука на биллиардной книжке $\Omega_n$, склеенной из $n$ экземпляров биллиарда в круге $D$ по перестановке $(1, 2, \dots, n)$, послойно гомеоморфна особенности фокус-фокус c $n$ положениями равновесия. Тем самым произвольная полулокальная особенность фокус-фокус реализована в классе биллиардных книжек с потенциалами. Теорема 2 доказана. $\square$.
	
	\section{Заключение}
	
	Тем самым, один из важных классов полулокальных особенностей --- невырожденные фокусные особенности --- удалось полностью реализовать в топологической категории интегрируемыми биллиардами. Развитие решенной нами задачи возможно, например, по следующим направлениям.
	
	Во-первых, особенности фокус-фокус в интегрируемых системах имеют нетривиальные гладкие и симплектические инварианты. Они устроены сложнее, чем топологические. Так гладкие инварианты фокусных особенностей с двумя положениями равновесия на одном слое описаны в работе А.В.Болсинова и А.М.Изосимова \cite{BolsIzosimov}. Симплектические инварианты фокусных особенностей полуторических интегрируемых систем ввел S.\ V\~u Ng\d oc в случае одной особой точки на слое \cite{Ngoc1} и (совместно c A.\ Pelayo) в общего случае для того же класса систем \cite{Ngoc2}. Вопрос о наличии их аналогов для кусочно-гладких биллиардов на столах-комплексах с потенциалом требует отдельного изучения.
	
	Другим обобщением может являться проблема реализации особенностей систем с большим чем 2 числом степеней свободы с помощью многомерных софокусных биллиардов с потенциалом. Такие топологические инварианты изучались, например, Н.Т.Зунгом \cite{14} и А.А.Ошемковым \cite{Osh10}, \cite{KozOsh20}, см. также \cite{2}).

Работа выполнена при поддержке гранта РНФ 17-11-01303 в МГУ имени М.\,В.~Ломоносова



%


\end{fulltext}

\end{document}